# GLOBAL COMMUNICATION ALGORITHMS FOR CAYLEY GRAPHS

*Vance Faber*
*Big Pine Key, Florida 33043*

**ABSTRACT**

We discuss several combinatorial problems that arise when one looks at computational algorithms for highly symmetric networks of processors. More specifically, we are interested in minimal times associated with four communication tasks (defined more precisely below): universal broadcast, every processor has a vector that it wishes to broadcast to all the others; universal accumulation, every processor wishes to receive the sum of all the vectors being sent to it by all the other processors; universal exchange, every processor wishes to exchange a vector with each other processor; and global summation, every processor wants the sum of the vectors in all the processors .

**§1. Introduction.** We discuss several combinatorial problems that arise when one looks at computational algorithms for highly symmetric networks of processors. More specifically, we are interested in minimal times associated with three communication tasks (defined more precisely below): universal broadcast, every processor has a vector that it wishes to broadcast to all the others; universal accumulation, every processor wishes to receive the sum of all the vectors being sent to it by all the other processors; and universal exchange, every processor wishes to exchange a vector with each other processor.

Algorithms for applying our results to hypercubes have been discussed in [1]. Several papers have looked at the problems of mapping communication algorithms onto multiprocessor machines. Reference [2] uses the notion of bisection to get lower bounds for communication tasks. Reference [3] discusses the minimum time to perform communication tasks that are permutations on networks whose graphs are stars. This paper is adapted from and corrects errors in an unpublished internal report from Los Alamos National Laboratory [10].

Our general model of a network is a directed graph with processors as vertices and the connections between them as edges. The case has been made elsewhere that a multiprocessor network should be homogeneous; that is, the network should appear the same from any processor. This means that the graph is vertex transitive. Sabiduissi [6] has shown that a graph is vertex transitive if and only if it is the Cayley coset graph of a group. If the coset is the identity group, the graph is called a Cayley graph. Our main intent is to create methods for scheduling highly symmetric global communication tasks on vertex symmetric networks. We discuss the notion of a regular order on a group of order $P$ with $d$ generators and we show that possessing a regular order is sufficient to ensure that if we have a network on the corresponding Cayley graph and we can use all the wires simultaneously on every time step, then the time for universal broadcast will be optimal, $\lceil (P-1)/d \rceil$. We show that the hypercube has these conditions.

We also discuss the difficulties that arise when wires exist in both directions between any two connected processors but both cannot be used simultaneously. In this case, proving any general theorems seems daunting. As an example, we analyze the case of the hypercube. We show that if all the wires on a $d$ dimensional hypercube can be used simultaneously on every time step but only in one direction then universal broadcast and universal accumulation take $\lceil 2(P-1)/d \rceil$ time steps while universal exchange takes $P$ time steps.

Finally, we can show that if the graph is *distance transitive* (see [12] and defined more precisely below), then global summation can be accomplished in the number of steps equal to the diameter of the graph. So in the *d* dimensional hypercube, for example, we can take each step to be the set of wires in a fixed coordinate direction. After *d* of these steps, all the processors have the sum of all the data.

**§2. Definitions and examples.**

*2.1. Definition* (*Cayley coset graph*). Let $\Gamma$ be a finite group, $H$ be a subgroup, and $\Delta$ be a set of distinct nonidentity coset representatives of $H$ in $\Gamma$ with the properties

(i) $\Delta \cup H$ generates $\Gamma$;
(ii) $H\Delta H = \Delta H$.

We define a Cayley coset graph $G = G(\Gamma, \Delta, H)$ with vertices $gH, g \in \Gamma$ (the left cosets of $H$ in $\Gamma$), and edges $(gH, g\delta H)$ for every $g \in \Gamma$ where $\delta \in \Delta$. (Note that the properties (i) and (ii) assure that $G$ is a connected directed graph with no loops or multiple edges which is regular with both out-degree and in-degree $d = |\Delta|$. All connected vertices are connected in both directions if and only if for each generator $\delta$ there is a generator $\delta'$ with $\delta^{-1}H = \delta'H$.)

The hypercube is a Cayley coset graph. The group is $\Gamma = Z_2 \times \ldots \times Z_2$, where $Z_2$ is the two-element group; $H$ is the identity group; $\Delta$ is the set of canonical generators $\{(1,0,\ldots,0),(0,1,\ldots,0),\ldots,(0,0,\ldots,1)\}$. Since the square of each generator is the identity, the hypercube has wires in both directions between any two connected processors. We shall say more about Cayley coset graphs in Section 4.

*2.2. Definition.* A *task graph T* is a sequence of directed edges $\{e_i \mid i \in I\}$ of a graph $G$ labeled by positive integers (called times) $t(e_i)$, satisfying

i) $t(e_i) < t(e_j)$ implies that either $e_i < e_j$ (see below) or $e_i, e_j$ are incomparable,

ii) $t(e_i) = t(e_j)$ implies that either $i = j$ or $e_i, e_j$ are incomparable.

*2.3. Definition.* We say that the edges $e_i, e_j$ satisfy $e_i < e_j$ if and only if there is a directed path in $T$ beginning with $e_i$ and ending with $e_j$. We say that $e_i$ and $e_j$ are *incomparable* if neither $e_i < e_j$ nor $e_j < e_i$.

*2.4. Definition.* The time $\tau(T)$ for a task graph is

$$\tau = \max_{i \in I} t(e_i).$$

Note that a task graph is a labeled digraph. As such, we can refer to it by its adjacency matrix or we might want to utilize one adjacency matrix for each of its $\tau$ time steps.

*2.5. Example (see Fig 1).* Four edges: $e_1 = (0,1)$, $e_2 = (1,2)$, $e_3 = (2,3)$, $e_4 = (3,0)$ with times $t(e_i) = i$. Even though $e_4 < e_1$, we do not require $t(e_4) < t(e_1)$.

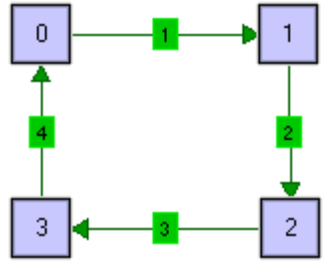

Figure 1

The four adjacency matrices of the time steps are

$$A_1 = \begin{bmatrix} 0 & 1 & 0 & 0 \\ 0 & 0 & 0 & 0 \\ 0 & 0 & 0 & 0 \\ 0 & 0 & 0 & 0 \end{bmatrix}, A_2 = \begin{bmatrix} 0 & 0 & 0 & 0 \\ 0 & 0 & 1 & 0 \\ 0 & 0 & 0 & 0 \\ 0 & 0 & 0 & 0 \end{bmatrix}, A_3 = \begin{bmatrix} 0 & 0 & 0 & 0 \\ 0 & 0 & 0 & 0 \\ 0 & 0 & 0 & 1 \\ 0 & 0 & 0 & 0 \end{bmatrix}, A_4 = \begin{bmatrix} 0 & 0 & 0 & 0 \\ 1 & 0 & 0 & 0 \\ 0 & 0 & 0 & 0 \\ 0 & 0 & 0 & 0 \end{bmatrix}.$$

*2.6. Example (see Fig 2).* Six edges: $e_1 = (0,1)$, $e_2 = (1,2)$, $e_3 = (1,3)$, $e_4 = (3,4)$, $e_5 = (2,4)$, $e_6 = (4,5)$ with times $t(e_1) = 1$, $t(e_2) = 2$, $t(e_3) = 2$, $t(e_4) = 4$, $t(e_5) = 3$, $t(e_6) = 5$. The relation $t(e_5) < t(e_4)$ is allowed since $e_5$ and $e_4$ are incomparable.

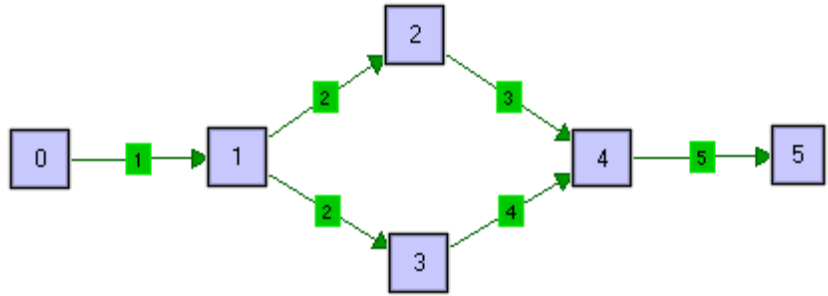

Figure 2

*2.7. Discussion.* How a communication task is associated with a task graph.

We use a simple model that assumes that when a word of data gets to a processor, it takes no time to store it, retrieve it, or operate on it. The only time that counts is the time to move one word of data across a wire (a fixed constant throughout the network). The task graph represents the time step in which a given piece of data transfers between two processors. In [9], more realistic models are discussed.

*2.8. Example.* The task graph for moving one word of data from a single processor 0 to all others in the cube $Q_3$ might be a directed tree that spans all processors. See Figure 3.

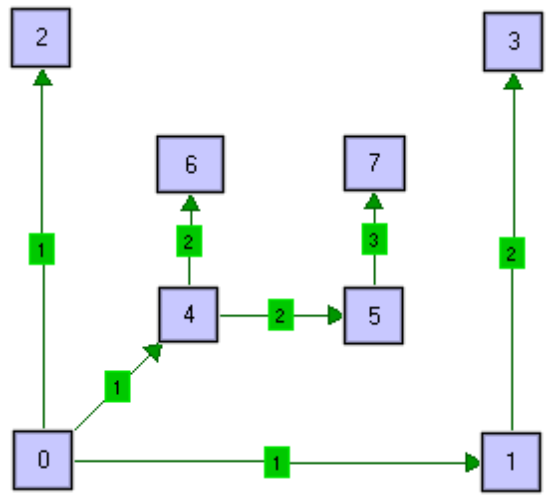

Figure 3.

*2.9. Example.* Reversing all arrows and times $(t' = \tau - t + 1)$ yields the task graph for accumulating one word from each of the processors and storing it in processor 0. See Figure 4.

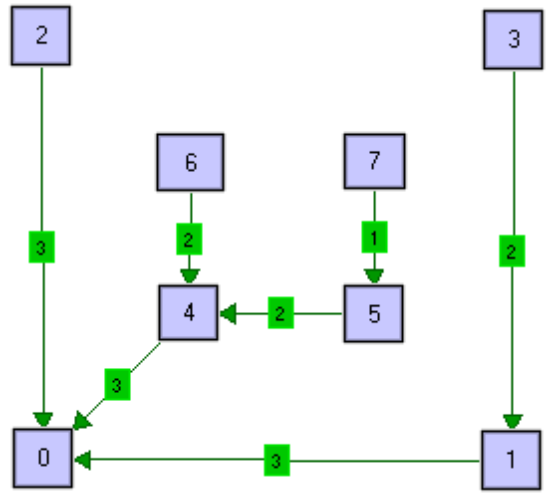

Figure 4.

*2.10. Theorem.* If $T$ is a task graph, the graph $T'$ formed by reversing all arrows and times is a task graph.

*Proof.* Clearly $e_i < e_j$ in $T$ if and only if $e_j < e_i$ in $T'$ and $t(e_i) < t(e_j)$ if and only if $t'(e_j) < t'(e_i)$.

*2.11. Theorem (Duality).* $T$ is a task graph for moving one word from one processor $P_1$ to all the others (a broadcast from $P_1$) if and only if $T'$ is a task graph for accumulating one word from each processor and storing it in $P_1$ (an accumulation to $P_1$).

*Proof.* $T$ will move one word from $P_1$ to all the others if and only if it has a directed path from $P_1$ to each of the others. This is equivalent to saying that $T'$ has a directed path from each processor to $P_1$. The condition in $T'$ that the sum of two (or more) incoming words cannot exit a processor before each summand has arrived is equivalent to the condition in $T$ that a word cannot be sent out to two (or more) processors before it arrives.

*2.12. Discussion.* Not every communication task can be represented as a task graph. Essentially, we have forced a task graph to pass only one word of information.. (It may be in slightly different form at different times.) For this reason, we need a more general definition. For example, in a hypercube if we were allowed to use the wires in both directions between processors simultaneously, the tasks in Figure 3 and Figure 4 could proceed simultaneously.

*2.13. Definition.* A *communication graph C* is a collection of task graphs having the property that no directed edge in *C* can have a given label more than once.

*2.14. Definition.* A *universal broadcast* is a communication graph consisting of a broadcast from each processor.

*2.15. Definition.* A *universal exchange* is a communication graph consisting of a directed path from each vertex to each other.

*2.16. Example.* A communication graph for universal broadcast in a 2-cube (see Figure 5):

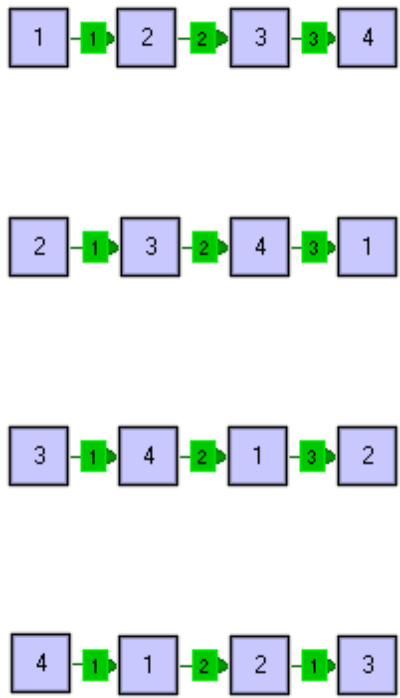

Figure 5.

*2.17. Definition.* The *time* for a communication graph is the maximum of the times for each of its tasks,

$$\tau(C) = \max_{T \in C} \tau(T).$$

*2.18. Discussion (one-way communication).* In our model of communications, we have specified that when there are edges in both directions between two processors, they both can be used on a single time step.. In some networks, both of these edges share the same physical wire and cannot be used in both directions at once. With this in mind, we make the following definition.

*2.19. Definition.* A *one-way* communication graph $C$ is a collection of task graphs having the property that edges between the same two vertices in $C$ cannot have the same label. To distinguish this model and the original one, we put the words "one-way" in front of the previous definitions.

*2.20. Definition.* One task graph for *global sum* is two spanning trees rooted in a single vertex. One tree has edges directed in to the vertex and the other out from the vertex. Times start on the leaves of the first tree and increase towards the root. The times then continue to increase on the second tree from the root to the leaves.

**§3. Lower bounds.**

*3.1. Theorem.* For a universal broadcast,

$$\tau \geq \frac{P(P-1)}{Q}$$

where $Q$ is the number of edges in $G$ which can be used simultaneously.

*Proof.* Each task graph has at least $P-1$ edges. There are $P$ task graphs, so we need $P(P-1)$ edges in $C$, at least. But $C$ can have at most $\tau Q$ edges (then every edge appears with all times). Thus $\tau Q \geq P(P-1)$.

*3.2. Theorem.* For a regular bi-directed graph of degree $d$, the time for a one-way universal broadcast is

$$\tau \geq \frac{2(P-1)}{d}.$$

*Proof.* Since edges can only be used in one direction at a time, $Pd = \sum_{i=1}^{P} d_i = 2Q$, so $Q = \frac{Pd}{2}$.

*3.3. Theorem.* For the *d*-dimensional hypercube $Q_d$, let $\tau$ be the time for a one-way universal broadcast from each point to all other points except the single point farthest away. Then

$$\tau \geq \frac{2(P-2)}{d}.$$

*Proof.* Each task graph has at least $P-2$ edges.

*3.4. Example.* If *d* is an odd prime, then (mod *d*)

$$2^d - 2 \equiv 2 \cdot 2^{d-1} - 2 \equiv 2 \cdot 1 - 2 \equiv 0$$

since for every prime $d$, $2^{d-1} \equiv 1$ (mod *d*). Thus it is possible that the minimum can be achieved.

*3.5. Theorem.*

(i) For a universal broadcast on a regular bi-directed graph,

$$\tau \geq \frac{P(P-1)}{2Q}$$

where *Q* is the number of bi-directed edges in *G*.

(ii) For a universal broadcast on a regular directed graph with both in-degree and out-degree *d*,

$$\tau \geq \frac{P-1}{d}.$$

*Proof.* Left to the reader.

**§4. Universal broadcast on Cayley coset graphs.** In this section, we describe a property possessed by many Cayley coset graphs that is a sufficient condition for optimal universal broadcast on graphs of groups.

*4.1. Definition.* If *G* is a directed graph and *v* is a vertex in *G*, then for any positive integer *r*, $S_r(v)$ is the set of all vertices *u* in *G* for which there is a directed path from *v* to *u* of length less than or equal to *r*.

*4.2. Remark.* Note that if $G = G(\Gamma, \Delta, H)$, then $gH \in S_r(H)$ if and only if *g* has a representation modulo *H* as a product $\pi$ of *r* or less elements from $\Delta$. Also, $gH \in S_r(H)$ if and only if *g* has the form $h\pi$ modulo *H* with $h \in H$.

*4.3. Definition.* We say that $G = G(\Gamma, \Delta, H)$ has a *regular order* if there is an indexing of $\Delta$, $\{\delta_1, \cdots, \delta_d\}$, and an indexing of the cosets of $H$ in $G$, $\{g_0 H, \cdots, g_{P-1} H\}$, so that

(i) $g_0 = e$;

(ii) if $g_j H \in S_r(H)$ and $i < j$, then

$$g_i H \in S_r(H);$$

(iii) fix $a \geq 0$ and $1 \leq b \leq d$ and consider a set of cosets of $H$ in $G$ of the form $g_{ad+c} H$, $1 \leq c \leq b$. Then there exists $\delta_{i_1}, \delta_{i_2}, \cdots, \delta_{i_b} \in \Delta$, all different, and such that for each $c$, $g_{ad+c} \delta_{i_c}^{-1} H = g_{s_c} H$ with $s_c \leq ad$ and if $g_{ad+c} H \in S_r \setminus S_{r-1}$ then $g_{s_c} H \in S_{r-1}(H)$.

*4.4. Explanation.* We can think of the elements of $\Delta$ as directions. What part (iii) of this definition says is that the vertices of $G$ can be broken up into blocks of $d$ in such a way that in each block, vertices are connected in distinct directions to vertices closer to the root vertex, $g_0 H = H$. Parts (i) and (ii) say that the blocks are ordered so that the vertices are listed in non-decreasing order with respect to distance from $g_0 H = H$. The main strength of regular ordering is that, as we shall now show, it is a sufficient condition in Cayley graphs for the existence of an optimal two-way universal broadcast.

*4.5. Theorem.* In a Cayley graph that has a regular order with degree $d$ and $P$ vertices, the optimal time for universal broadcast is

$$\left\lceil \frac{P-1}{d} \right\rceil.$$

*Proof.* The lower bound is the content of Theorem 3.5(ii). Using the labeling in Definition 4.3(iii), consider the tree $T_{g_0}$ whose edges are $(g_{ad+r} \delta_{i_r}^{-1}, g_{ad+r}) = (g_{s_r}, g_{ad+r})$ for all $a$ and $r$ so that $1 \leq r \leq d$ and $ad + r \leq P - 1$. We assign the time $a + 1$ to all the edges $(g_{s_r}, g_{ad+r})$. In the tree $T_{g_i}$ with edges $(g_i g_{s_r}, g_i g_{s_r} \delta_{i_r}) = (g_i g_{s_r}, g_i g_{ad+r})$, we also assign the time $a + 1$ to the edges $(g_i g_{s_r}, g_i g_{ad+r})$. We claim that $C = \bigcup_{i=0}^{n-1} T_{g_i}$ is a communication graph for universal broadcast. Parts (i), (ii) and (iii) of Definition 4.3 guarantee that $T_{g_i}$ is a task graph. Suppose two of the edges labeled $a + 1$ in $C$, say $(g_i g_{s_r}, g_i g_{s_r} \delta_{i_r}) = (g_i g_{s_r}, g_i g_{ad+r})$ and $(g_j g_{s_t}, g_j g_{s_t} \delta_{i_t}) = (g_j g_{s_t}, g_j g_{ad+t})$, are identical. Two edges in a Cayley graph cannot be the same unless they emanate from the same group element and arise from the same generator. Thus $\delta_{i_r} = \delta_{i_t}$ and $g_i g_{ad+r} = g_j g_{ad+t} = g_j g_{ad+r}$, so $g_i = g_j$; that is, the two edges are not in different task graphs. This contradiction proves the theorem.

*4.6. Definition.* Good communication networks often enjoy additional symmetries besides vertex transitivity. A graph is distance transitive if for every two pairs of vertices $(v, w)$ and $(x, y)$ with common distances $d(x, y) = d(v, w)$ there exists an automorphism of the graph $\theta$ with $\theta(v) = x$ and $\theta(w) = y$.

*4.7. Remarks.* The power of a Cayley graph for universal communication lies in the fact that if $T_{g_0}$ is a task graph which is a directed tree rooted at $g_0$, then the tree $T_{g_i}$ is a task graph rooted at $g_i$. We call $T_{g_0}$ a *template*. We shall show in Section 7 that the hypercube can be regularly ordered. We conjecture that if $(\Gamma_1, \Delta_1, E)$ and $(\Gamma_2, \Delta_2, E)$ can be regularly ordered, so can

$$(\Gamma_1 \times \Gamma_2, \Delta_1 \times \{e\} \cup \{e\} \times \Delta_2, E).$$

Do all Cayley coset graphs admit optimal universal broadcast? If a vertex symmetric graph can be regularly ordered, then clearly $D \leq \left\lceil \frac{P-1}{d} \right\rceil$. Example 4.7 below shows that is not the case. Can all vertex transitive graphs with $D \leq \left\lceil \frac{P-1}{d} \right\rceil$ be regularly ordered? Perhaps the time for universal broadcast on a vertex symmetric graph is the maximum of $\left\lceil \frac{P-1}{d} \right\rceil$ and $D$. We have not been successful dealing with Cayley coset graphs whether they can be regularly ordered or not. In [10], we specified an additional property of Cayley coset graphs that we claimed would be sufficient for optimal universal broadcast but Ken Blaha [11] showed that the Petersen Graph had a regular order with this additional property but the task graphs produced by the algorithm had identical edges assigned the same time. It is easy to see that universal broadcast on the Petersen Graph can be done in the time predicted by Theorem 4.5 (three time steps) so perhaps there is still some way to apply this theory for Cayley coset graphs.

*4.7. Example.* Take the abelian group $Z_2 \times Z_8$ with the 5 generators $(1,0)$, $(-1,0)$, $(0,1)$, $(1,1)$, $(-1,1)$. The distance between $(0,0)$ and $(4,0)$ is 4 while $\frac{P-1}{d} = 3$.

**§5. Global sum on Cayley coset graphs.** Let $A$ be the adjacency matrix of a graph. Suppose we use the entire graph as a single time step $t$ of a communication graph. If $x$ is a vector such that processor $i$ has a value $x_i$ to communicate at time $t$, then $(A - cI)x$ is the vector obtained after each processor adds the values it gets on time $t$ and subtracts $c$ times its own value. The global sum is the inner product $(x, \vec{1})$ where $\vec{1}$ is the vector with all entries equal to one.

*5.1. Theorem.* Suppose that the regular graph $A$ has exactly $D+1$ distinct eigenvalues. Then the time for global sum is at most $D$.

*Proof.* We choose time steps of the form $(A - \lambda_t I)$ where the $\lambda_t$ are the eigenvalues other than $d$ (which corresponds to the constant eigenvector $e = \vec{1}/\sqrt{n}$). After the $k^{th}$ time step, the $i^{th}$ processor has the $i^{th}$ entry of the vector $\prod_{t=1}^{k}(A - \lambda_t I)x$. Now if we expand $x$ into its parts along $e$ and orthogonal to $e$, we have $x = y + (x,1)\vec{1}/n$ where the sum of the entries of $y$ is zero. Furthermore $\prod_{t=1}^{D}(A - \lambda_t I)y = 0$, so after $D$ time steps, the entries are all the same value

$$\mu = \frac{(x,1)}{n} \prod_{t=1}^{D}(d - \lambda_t).$$

We recover the global sum at each processor by scaling $\mu$ by the predetermined factor.

*5.2. Corollary.* If all the edges in a distance transitive graph can be used simultaneously, then global sum takes $D$ time steps where $D$ is the diameter (the greatest distance between two vertices).

*Proof.* A distance transitive graph has exactly $D+1$ eigenvalues (see [12, page 113]).

*5.3. Example.* The hypercube has exactly $d+1$ distinct eigenvalues diameter $d$ so global sum can be accomplished in $d$ time steps.

**§6. The hypercube.** In this section, we we state theorems specifically for the hypercube and one-way communication. The proofs are carried out in the next section. Since we shall show that the hypercube can be regularly ordered, universal broadcast can be accomplished in

$$\left\lceil \frac{P-1}{d} \right\rceil = \left\lceil \frac{2^d - 1}{d} \right\rceil$$

time steps by Theorem 4.5. Since there are edges in both directions between every pair of connected vertices in the hypercube, we can split the edges into two sets with oppositely directed edges in different sets. For each time $t$ in the universal broadcast, we assign a new time $2t - 1$ to an edge if it was in the first set and $2t$ if it was in the second set. This produces a one-way universal broadcast with time

$$2\left\lceil \frac{P-1}{d} \right\rceil = 2\left\lceil \frac{2^d - 1}{d} \right\rceil.$$

This differs from the optimum by at most one depending on the remainder inside the floor function. It is quite complicated to deal with this small difference. At any rate, we do not know how to approach universal exchange in a general Cayley graph so specialized proofs are required for the hypercube in that case. Here are the theorems that we shall prove in the next section.

*6.1. Theorem.* For the $d$-cube, let $\tau$ be the time for a one-way universal broadcast to all but the point farthest away. Then

$$\tau \geq \frac{2(P-2)}{d}.$$

*6.2. Example.* If $d$ is an odd prime, then $(\bmod\, d)$

$$2^d - 2 \equiv 2 \cdot 2^{d-1} - 2 \equiv 2 \cdot 1 - 2 \equiv 0$$

since for every prime $d$, $2^{d-1} \equiv 1 \,(\bmod\, d)$. Thus it is possible that the minimum can be achieved.

*6.3. Theorem.* The time for a one-way universal broadcast in a $d$-cube is

$$\left\lceil 2\frac{2^d - 1}{d} \right\rceil.$$

In fact, the time for a broadcast from each vertex to all those within distance $l$ is

$$\left\lceil \left(2\sum_{i=1}^{l}\binom{d}{i} - 2\right)/d \right\rceil.$$

In addition, on odd time steps, the edges can be directed arbitrarily, with the opposite direction used on the next time step.

*6.4. Remark.* Reference [7] has basically the same algorithm, but seems to avoid all the complications by restricting attention only to cubes whose dimensions are prime numbers and then the fraction is an integer.

*6.5. Theorem.* The time for one-way universal exchange in a $d$-cube is $2^d$. In fact, the time to communicate from each vertex simultaneously to the vertices at distance $d-1$ and $d$ is $2d$, while the time to communicate from each vertex to the vertices at distance $s \leq d-1$ is $2\binom{d-1}{s-1}$. In addition, on odd time steps the edges can be directed arbitrarily, with the opposite direction used on the next time step.

*6.6. Definition.* The communication graph consisting of a broadcast from each processor to those within distance $l$ is called *universal broadcast to within distance $l$*. The communication graph consisting of a directed path from each vertex to those within distance $l$ is called *universal exchange to within distance $l$*. The communication graph formed from the universal broadcast to distance $l$ by the process in Theorem 2.10 is called *universal accumulation from within distance $l$*.

*6.7. Remark.* We shall find optimal one-way times for all of these communication tasks on the hypercube. Lower bounds are always found from the simple counting arguments as in Theorems 3.1 and 3.2; these are omitted from now on.

## §7. Hypercube proofs.

*7.1. Lemma.* Let $D$ be a set with $d$ elements. Let $(e_1, e_2, \cdots, e_d)$ be a list of the elements of $D$ in some order. Let $s$ be a nonzero integer less than or equal to $d$. The $s$-element subsets of $D$ can be ordered $(S_i \mid 1 \leq i \leq \binom{d}{s})$ so that $e_j \in S_i$ if $i \equiv j \pmod{d}$.

*Proof.* Let $P$ be the permutation which cyclically permutes $(e_1, e_2, \cdots, e_d)$; that is, $Pe_i = e_{i+1}$ with the indices read modulo $d$. We can partition the set $D^{[s]}$ of $s$-element subsets into equivalence classes by the equivalence relation $S \sim T$ if and only if $P^k S = T$ for some $k$. Note that if $e_k \in S$, then $e_{k+i} \in P^i S$ for each $i$. Thus, if an equivalence class has $d_1$ elements, $S, PS, \cdots, P^{d_1-1}S$, some element contains a given $e_k$ say $S_1$ and it is easy to order the elements of the class so that $e_{k+i} \in P^i S_1 = S_{i+1}$, $i = 0, 1, \cdots, d_1 - 1$. To order all of $D^{[s]}$ correctly, we proceed through the classes ordering each one as above with $k$ taken in each instance to be the next (modulo $d$) after the last index used in the previous step.

*7.2. Lemma.* Let $D$ be a set with $d$ elements Let $(e_1, e_2, \cdots, e_d)$ be a list of the elements of $D$ in some order. The nonempty subsets of $D$ can be ordered, $(S_i \mid 1 \leq i \leq 2^{d-1})$, so that $e_j \in S_i$ if $i \equiv j \pmod{d}$ and if $i \leq j$ then $|S_i| \leq |S_j|$.

*Proof.* Obvious.

*7.3. Definition.* We say that the edges $(u, v)$ and $(u', v')$ in $Q_d$ are parallel if $u + v = u' + v'$. We call $u + v = e_i$, for some $i$, the direction of $(u, v)$.

*7.4. Lemma.* Let $S = \{(u_i, v_i) : i = 1, 2, \cdots s\}$ be a set of edges in $Q_d$ having the property that no two are parallel. Assume $e_i = u_i + v_i$. (One can think of the vertices of the cube as subsets of $(e_1, e_2, \cdots, e_d)$; then the plus sign denotes symmetric difference of sets. Alternatively, one can think of the vertices as vectors in $Z_2^d$; then the plus sign denotes vector addition and $e_i$ is the canonical unit vector in the $i^{th}$ direction. We use these concepts interchangeably.) Suppose $X$ is a set of vertices in $Q_d$ satisfying $x, y \in X$ implies that $x + y \neq e_i$, $i = 1, 2, \cdots s$. Let $S + X = \{(u_i + x, v_i + x) : x \in X, 1 \leq i \leq s\}$. Then an edge $(z, z + e_k)$ appears in the sequence $S + X$ at most once. In fact, $(z, z + e_k)$ appears in $S + X$ if and only if $1 \leq k \leq s$ and either $z + u_k \in X$ or $z + v_k \in X$.

*Proof.* If $(z, z + e_k)$ appears in $S + X$, then either

1) $(z, z + e_k) = (x + u_i, x + v_i)$

or

2) $(z + e_k, z) = (x + u_i, x + v_i)$.

If 1) holds, then

$$\begin{cases} z = x + u_i \\ z + e_k = x + v_i \end{cases}$$

so $e_i = u_i + v_i = e_k$.

If 2) holds, then

$$\begin{cases} z + e_k = x + u_i \\ z = x + v_i \end{cases}$$

so $e_k = u_i + v_i = e_i$.

Thus, we have two possible solutions, both of which require $1 \leq k \leq s$. In the first $x_1 = z + u_i \in X$ and in the second $x_2 = z + v_i \in X$. If both of these could hold, we would have $x_1 + x_2 = u_i + v_i = e_i$, violating our assumption on $X$.

*7.5. Corollary.* Let $S = \{(u_i, v_i) : i = 1, 2, \cdots s\}$ be a set of edges in $Q_d$ having the property that no two are parallel. Then the edge $(z, z + e_k)$ which appears in $S + Q_d$ appears exactly twice.

*Proof.* Let $E, O$ be the even and odd vertices, respectively, where the even vertices correspond to sets with an even number of elements and the odd vertices correspond to sets with and odd number of elements. Then by the Lemma with $X = E$ (or $O$), $(z, z+e_k)$ appears in $E$ (or $O$) if and only if $1 \le k \le s$ and either $z+u_k \in X$ or $z+v_k \in X$. Thus $(z, z+e_k)$ appears in $S+E$ if and only if it appears in $S+O$.

*Proof of Theorem 6.3.* We order the nonzero elements of $Q_d$ as in Lemma 7.2. Consider only those subsets $S_i, i \in I$, with at most $l$ elements. Let $v_i = \sum_{k \in S_i} e_k$. We first form the directed graph $T_0$ with edges $(u_i, v_i)$ with $u_i = v_i + e_{i(\bmod d)}$. This is the task graph for broadcasting to all vertices of distance $l$ from vertex 0 of the hypercube when labeled as follows. We can write $i = ad + b$ with $1 \le b \le d$ and $0 \le a$. Then the labeling $t$ is given by $t(u_i, v_i) = 2a+1$. We shall show below that "usually" $T_0$ labeled in this way is a task graph. For the moment, let us assume that $T_0$ is a task graph. To form the communication graph for universal broadcast to distance $l$, we form a directed tree $T_x$ for each nonzero $x \in Q_d$ whose directed edges are $(x+u_i, x+v_i)$. We label $T_x$ as follows:

$$t(u_i + x, v_i + x) = \begin{cases} 2a+1 & \text{if } x \text{ is even} \\ 2(a+1) & \text{if } x \text{ is odd} \end{cases}.$$

It is easily seen that if $T_0$ is a task graph, then so is $T_x$. We claim that the collection of graphs $C = \{T_x, x \in Q_d\}$ is a communication graph for universal broadcast to distance $l$. To see this, let $S_a$ be the set of edges in $T_0$ labeled with time $2a+1$. If we apply Lemma 7.4, we find that each edge in $S_a + E$ appears at most once and is labeled $2a+1$, while each edge in $S_a + O$ appears at most once and is labeled $2(a+1)$. Clearly $S_a + E = \{e \in T_x : x \text{ even}\}$ and $S_a + O = \{e \in T_x : x \text{ odd}\}$. Thus $C$ is a communication graph. Since $T_0$ connects all vertices of distance $l$ from vertex 0 to vertex 0 by a spanning tree and since $u+x$ is the same distance from $x$ as $u$ is from 0, $C$ is a communication graph for universal broadcast to distance $l$.

Let $N_l$ be the number of nonzero vertices in $Q_d$ at distance at most $l$ from vertex 0. We see from the labeling of $T_x$ that the time $\tau(C)$ for $C$ is at most $\alpha = 2 \left\lceil \frac{N_l}{d} \right\rceil$. This is not quite the value, $\tau(C) = \left\lceil 2 \frac{N_l}{d} \right\rceil$, given in the statement of the theorem. We shall now discuss this problem. Let $N_l \equiv r \pmod{d}$. The number $\alpha$ is one larger than $\left\lceil 2 \frac{N_l}{d} \right\rceil$ just in case $1 \le r \le d/2$. The problem arises because there are only $r$ edges labeled $\alpha - 1$ in $T_0$. These $r$ edges lead to $r$ edges labeled $\alpha - 1$ in $T_x$ for $x$ even. Thus only $r2^d/2$ edges in $Q_d$ are labeled with time $\alpha - 1$. This same number of edges is labeled with time $\alpha$; namely, those edges derived from the edges labeled $\alpha - 1$ in $T_0$ by adding $x$ when $x$ is odd. We would like to label these latter edges with $\alpha - 1$ also; however, they are exactly the same edges since they use the same $r$ directions. To fix this, we need to have the last $r$ subsets, $R_1, R_2, \cdots, R_r$, given by the ordering in Lemma 7.2, have the property that not only is $e_i \in R_i, i = 1, 2, \cdots, r$ but $e_{i+r} \in R_i, i = 1, 2, \cdots, r$. We shall leave the proof of the possibility of an ordering of this type to a technical lemma (Lemma 7.6) given below. Suppose the ordering we are using to construct $T_0$ has this property. Let $v_i = \sum_{j \in R_i} e_j$, $i = 1, 2, \cdots, r$. We need to have our ordering possess one additional property. If we form $T_0'$ from $T_0$ by removing the edges $(v_i + e_i, v_i)$ and replacing them with $(v_i + e_{i+r}, v_i)$, we need to have $T_0'$ be a task graph. This is related to the problem of having $T_0$ itself be a task graph, so we discuss both of these problems in a moment. Now we form $T_x'$ for $x$ odd by removing the edges $(v_i + e_i + x, v_i + x)$ from $T_x$ and replacing them with $(v_i + e_{i+r} + x, v_i + x)$. It is clear that we now have a communication graph for universal broadcast to distance $l$ which has time $\left\lceil 2 \frac{N_l}{d} \right\rceil$.

Now we must discuss the problem of $T_0$ (and $T_0'$) being task graphs. The edges $(v_i + e_i, v_i)$ in $T_0$ have the time $2a + 1$ just in case $i = ad + b$ with $1 \leq b \leq d$. For fixed $a$, consider the $i$ of this form. $T_0$ will be a task graph as long as $u_i = v_i + e_i$ is not one of the $v_i$, otherwise $u_i$ will not acquire the necessary information in time to transmit it to $v_i$; that is, the ordering requirement in the definition of task graph will not be satisfied. Now, in general, each $v_i = \sum_{j \in S_i} e_j$ is associated with a set $S_i$ which has a fixed constant number of elements for all $b$, $1 \leq b \leq d$. In this case, $u_i$ cannot be one of the $v_i$ since it is associated with a set with fewer elements. In rare cases, however, the number of elements in the $S_i$ may vary by 1. We need an additional constraint on the ordering constructed in Lemma 7.2. Not only do we need to have $|S_i| \leq |S_j|$ when $i \leq j$, but if $j = ad + b$ with $1 \leq b \leq d$, we cannot have $S_j \setminus e_b = S_i$ with $i + ad + b'$. A similar property must hold for the special ordering used to construct $T_0'$. For each $j$ with $1 \leq j \leq r$, we cannot have either $R_j \setminus e_j$ or $R_j \setminus e_{j+r}$, be $R_i$ for some $i$ with $1 \leq i \leq r$. We leave the proof of the possibility of ordering of this type to a technical lemma (Lemma 7.7) given below.

Finally, note that the labeling $t$ has the flexibility that initially we could have either $t(u_i, v_i) = 2a + 1$ or $t(v_i, u_i) = 2a + 1$ as long as the opposite direction is labeled $2(a + 1)$. This proves the theorem.

**7.6. Lemma.** Suppose $N_l \equiv r \pmod{d}$ with $1 \leq r \leq d/2$. The non-empty subsets of $D$ with at most $l$ elements can be ordered $(S_i \mid 1 \leq i \leq N_l)$ so that $e_j \in S_i$ if $i \equiv j \pmod{d}$, if $i \leq j$ then $|S_i| \leq |S_j|$, and if $R_1, R_2, \cdots, R_r$ denotes the last $r$ subsets in the ordering, then $e_{j+r} \in R_j$, $j = 1, 2, \cdots, r$.

*Proof.* Clearly, if $r \neq 0$, then since $\binom{d}{1} = d$, we must have $l \geq 2$. Suppose $l < d$, then all of the $R_j$ have $l$ elements. We can easily construct a set $R_1$ which has the following properties: $e_1 \in R_1$, $e_{1+r} \in R_1$, and the equivalence class of $R_1$ under the action of the permutation $P$ (see the proof of Lemma 7.1) has at least $r$ distinct members. We can accomplish this last part by forcing $R_1$ to have its remaining $l-2$ elements so that they and $e_{1+r}$ have consecutive indices. Now let $R_k = P^{k-1}R_1$. We can now arrange the construction in the proof of Lemma 7.1 so that the last $r$ sets are $R_1, R_2, \cdots, R_r$. Finally, suppose that $l = d$ and $r > 1$. We construct $R_1, R_2, \cdots, R_{r-1}$ as above and then let $R_r = D$).

**7.7. Lemma..** The non-empty subsets of $D$ with at most $l$ elements can be ordered $(S_i \mid 1 \leq i \leq N_l)$ so that $e_j \in S_i$ if $i \equiv j \pmod{d}$, if $i \leq j$ then $|S_i| \leq |S_j|$, and if $j = ad + b$ with $1 \leq b \leq d$, then we do not have $S_j \setminus e_b = S_i$ with $i + ad + b'$. In addition, if $N_l \equiv r \pmod{d}$ with $1 \leq r \leq d/2$ and $R_1, R_2, \cdots, R_r$ denote the last $r$ subsets in the ordering, then $e_{j+r} \in R_j$, $j = 1, 2, \cdots, r$ and $R_j \setminus e_{j+r}$ is not equal to any $R_i$ with $1 \leq i \leq r$.

*Proof.* Suppose $l \neq d$. Again, we use the construction in the proof of Lemma 7.1 in reverse. Suppose $N_{l-j+1} = a_j d + r_j$ with $1 \leq r_j \leq d$. We order the $l$-element subsets as in Lemma 7.1 by using the list $(e_{r_1+1}, e_{r_1+2}, \cdots, e_d, e_1, \cdots, e_{r_1})$. If $S_1^1$ is the first set in this ordering, then $S_1^1 \setminus e_{r_1+1}$ has $l-1$ elements. We order the $(l-1)$-element subsets by using the list $(e_{r_2+1}, e_{r_2+2}, \cdots, e_d, e_1, \cdots, e_{r_2})$ and using the equivalence class of $S_1^1 \setminus e_{r_1+1}$ (under the shift operator $P$) at the very beginning of the ordering. Continue this process backwards until the ordering is completed. If there is more than one equivalence class of sets for a given $l - j + 1$, then $S_1^j \setminus e_{r_j+1}$ (where $S_1^j$ is the first set in the ordering of the $(l-j+1)$-element subsets) will not have the same time label as $S_1^j$ since between $S_1^j \setminus e_{r_j+1}$ and $S_1^j$ there are at least $d$ subsets. There are only two cases where there can be only one equivalence class. The first case is $l - j + 1 = 1$; however, here $r_1 = d$ and so clearly we have correctly ordered the one element subsets as $(e_1, e_2, \cdots, e_d)$. The second case is $l - j = d - 1$, but this can only happen when $l = d - 1$ and $j = 1$ and we have forced $S_1^j \setminus e_{r_j+1}$ (and its shifts) to have $l - 1$ elements. The construction in Lemma 7.6 guarantees the special requirements when $r_1 \leq d/2$. Finally, let us consider $l = d$. We wish to augment the construction used in the proof of Lemma 7.6. Again, we use $R_1, R_2, \cdots, R_{r_1} = D$ to denote the last $r_1$ subsets in the ordering. In addition, we denote by $T_{r_1}, T_{r_1+1}, \cdots, T_d$ the remaining $(d-1)$-element subsets which appear in the ordering directly before $R_1$. To satisfy the conditions of the Lemma, we need $e_{r_1} \in T_{r_1}, e_{r_1+1} \in T_{r_1+1}, \cdots, e_d \in T_d$ and both $e_1, e_{r_1} \in R_1$, both $e_2, e_{r_1} \in R_2$, ..., and both $e_{r_1-1}, e_{r_1} \in R_{r_1-1}$. The membership of $e_{r_1}$ in $R_1, R_2, \cdots, R_{r_1-1}$ is required so that $R_{r_1} \setminus e_{r_1} = D \setminus e_{r_1}$ is not one of the $R_1, R_2, \cdots, R_{r_1-1}$. In addition, if $1 \leq r \leq d/2$, then we also need $e_{1+r_1} \in R_1, e_{2+r_1} \in R_2, \cdots, e_{2r_1-1} \in R_{r_1-1}$. If $r_1 \neq d$, these conditions are easily satisfied, but if $r_1 = d$, they cannot be satisfied since they say that we must have $e_d \in T_d, e_d \in R_1, \cdots, e_d \in R_{d-1}$. In other words, when $r_1 = d$ in order to satisfy the conditions of the lemma, $e_d$ is required to be a member of all $(d-1)$-element subsets, clearly an impossibility. Now $r_1 = d$ means that $d$ divides $2^d - 1$. We shall show in the next lemma that this never happens.

*7.8. Lemma.* For all $d > 1$, $d$ cannot divide $2^d - 1$.

*Proof.* Suppose $d$ divides $2^d - 1$. Then clearly $d$ is odd. We use Fermat's theorem. Since $2^p - 2 \equiv 0 \pmod{p}$ for every prime $p$, $d$ is not a prime. Let $p_1, p_2, \cdots, p_n$ with $p_1 \leq p_2 \leq \cdots \leq p_n$ be the prime factorization of $d$. Then by assumption, $(2^{d/p_1})^{p_1} - 1 \equiv 0 \pmod{p_1}$. Since $(2^{d/p_1})^{p_1} - 2^{d/p_1} \equiv 0 \pmod{p_1}$ subtracting the two congruences we get $2^{d/p_1} \equiv 1 \pmod{p_1}$. Let $l$ be the order of 2 in the multiplicative group of $Z_{p_1}$. Then $l$ divides $p_1 - 1$. Hence $l \geq p_2$ and $l \leq p_1 - 1$, which is impossible.

Now we move on to the proof of theorem 6.5.

*7.9. Definitions.* We say that a subset $S$ of $D$ is regular if the equivalence class of $S$ under the cyclic permutation $P$ (see the proof of Lemma 7.1) has $d$ members. Otherwise, we call $S$ special. Let $b$ be the smallest non-negative integer so that $P^b S = S$. We call $b$ the block size of $S$. We identify the set $S$ with the {0,1} word which represents its characteristic function. The action of $P$ on this characteristic function is to shift all the bits to the right with the last bit wrapping around to the beginning of the word. Thus we see that if $S$ has block size $b$ then the word associated with $S$ is composed of $n = d/b$ blocks of identical strings of bits of length $b$. If $e_i \in S$, we call $S \setminus e_i$ an antecedent of $S$.

*7.10. Lemma.* If $S$ has more than one element it has a regular antecedent.

*Proof.* By applying $P$ to $S$, we may assume without loss of generality that $e_1 \in S$. We proceed by contradiction; assume every antecedent of $S$ is special. Let $i_2$ be the first integer greater than 1 such that $e_{i_2} \in S$. Then $S^* = S \setminus e_{i_2}$ is special with block size $b < d$. We may assume that $b$ is maximum among all antecedents of $S$. Let $w$ be the first block in the word that represents $S^*$, that is, the word composed of the first $b$ bits of $S^*$. Then $S^*$ is represented by the word $ww\cdots w$ with $w$ repeated $n = d/b$ times. Since the first bit of $S^*$ is 1 and the second non-zero bit of $S^*$ is beyond coordinate $i_2$, $w$ has at least $i_2$ bits and the $i_2^{th}$ bit is 0. Let $w_*$ be the word formed from $w$ by switching the first bit to a 0. Let $w^*$ be the word formed from $w$ by switching the $i_2^{th}$ bit to a 1. Then $\lambda = w^* w \cdots w_*$ is an antecedent of $S$ and so must have a block size $b_*$ with $b \leq b_* < d$. The junction between the first $b_*$ bits of $\lambda$ and the remaining bits falls somewhere after $w^*$ and before $w_*$. Thus we can write $w = w_1 w_2$, $\lambda = w^*(w_1 w_2)(w_1 w_2) \cdots (w_1 w_2) w_*$ and the first block $\phi$ of $\lambda$ is $w^*(w_1 w_2) \cdots (w_1 w_2) w_1$.

Then the last block $\psi$ of $\lambda$ is $\mu(w_1 w_2) \cdots (w_1 w_2) w_*$ where $w = \nu\mu$ and $\mu$ has the same length as $w_1$. If the interior pairs, $w_1 w_2$, are actually present in $\phi$, then we reach a contradiction by noting that the end of $\phi$, $w_1 w_2$, has one more 1 than the corresponding end of $\psi$, $w^*$. Otherwise $\phi = w^* w_1$ and $\psi = \mu w_*$. Again, $w_*$ has at least one less 1 (and possibly 2 less) than the corresponding end of $\phi$.

*7.11. Lemma.* There is a one-to-one mapping $\Theta$ from the set of special equivalence classes of $D^{[s]}$ under $P$ with $s < d$ into the set of subsets of regular equivalence classes such that $\Theta$ has the properties

(i) $\Theta[S]$ has $n-1$ elements, when $S$ has block size $b$ and $b = d/b$ blocks;
(ii) $\Theta[S_1] \cap \Theta[S_2] = \emptyset$ when $[S_1] \neq [S_2]$.

*Proof.* In other words, we are assigning to each special class a set of $n-1$ unique regular classes. The proof is just a matter of showing that there are enough regular classes compared to the number of special classes. First we estimate the number $p_n$ of sets with block size $b = d/n$. An overestimate for $p_n$ is given by having $n$ identical blocks of length $b$ with $t = s/n$ elements in each block. This is an overestimate because some of these might have block size smaller than $b$. The number of such sets is $\rho_n = \binom{b}{t}$. An underestimate for $p_1$ is then given by subtracting this number from the total number of $s$ element sets. Thus

$$p_1 \geq \rho_1 - \sum_{n \neq 1} \rho_n.$$

To estimate the number of equivalence classes of sets with $n$ blocks, $\varepsilon_n$, we note that such an equivalence class has $d/n$ members. Thus $\varepsilon_n = \dfrac{n}{d} p_n$. We are trying to show that

$$\varepsilon_1 \geq \sum_{n \neq 1} n \varepsilon_n .$$

Since $\varepsilon_n = 0$ unless $n$ divides both $d$ and $s$, it will suffice to show

$$\frac{1}{d}\left(p_1 - \sum_{n|(d,s)} p_n\right) \geq \sum_{n|(d,s)} \frac{n}{d}(n-1)p_n .$$

Rearranging gives

$$p_1 \geq \sum_{n|(d,s)} (n^2 - n - 1)p_n .$$

This is what we shall show.

First, notice that if we take just those sets of $s$-element subsets of $D$ formed by taking all possible arrangements of $n$ blocks of length $d/n$ with $s/n$ elements in each block, their number is $p_n^n$. Thus

$$p_1 \geq p_n^n .$$

Let us see if we can establish

$$p_1 \geq \sum_{n|(d,s)} (n^2 - n - 1)p_1^{1/n} . \tag{7.12}$$

To do this, we examine those cases when $n(n-1)p_1^{1/n} > 3p_1^{1/2}$. By elementary calculations, this holds only when $n = s = 3$, $d = 6$ or $9$; $n = s = 4$, $d = 8$; $n = s = 5$, $d = 10$; $n = s = 6$, $d = 12$. In these cases, Eq. (7.12) holds by inspection. Otherwise, we have

$$\frac{9}{4} s^2 \leq \binom{d}{s} \text{ for } d \geq 2s, \text{ so}$$

$$\sum_{n|(d,s)} (n^2 - n - 1) p_1^{1/n} \leq \frac{3}{2} s p_1^{1/2} \leq p_1 ,$$

as claimed.

*Proof of Theorem 6.5.* For each equivalence class of $D^{[s]}$ under $P$ choose a representative $S$. Let $b$ be the block size of $S$ and $n = d/b$. For $1 \leq i \leq n$, let

$$S_i^0 = S \cap \{e_j \mid (i-1)b + 1 \leq j \leq ib\}.$$

Note that $P^b S_i^0 = S_{i+1}^0$ with the indices read modulo $n$. Thus all the sets $S_i^0$ have the same cardinality and we can let $t = |S_i^0|$. Choose a sequence of sets $S_i^j$, $1 \leq j \leq t$, with $S_i^j$ a regular antecedent of $S_i^{j-1}$ and $S_i^t = \emptyset$. We from a directed graph $(T_S)_0$ whose edges have the form

$$\left( P^k \left( \bigcup_{m=1}^n S_i^{j_m} \right), P^k \left( \bigcup_{m=1}^n S_i^{k_m} \right) \right),$$

where $1 \leq k \leq b$ and $(k_m)_{m=1}^n$ immediately precedes $(j_m)_{m=1}^n$ in the *barber pole* ordering. We mean by this that each $(k_m)_{m=1}^n$ or $(j_m)_{m=1}^n$ is a vector of the form

$$j_m = \begin{cases} j & 1 \leq m \leq l \\ j+1 & l+1 \leq m \leq n \end{cases}$$

for some $j$ with $1 \leq j \leq t-1$. We say that $(k_m)_{m=1}^n$ immediately precedes $(j_m)_{m=1}^n$ in this ordering if $\sum_{m=1}^n k_m + 1 = \sum_{m=1}^n j_m$. For example, with $t = 2$, $n = 5$, we have 00000, 00001, 00011, 00111, 01111, 11111, 11112, 11122, 11222, 12222,, 22222. The graph $(T_S)_0$ has $dt = bnt$ edges which appear in groups of $b$; $b$ edges connecting the members of the equivalence class of $S$ to sets with $s-1$ elements (these edges have directions $e_1, e_2, \cdots, e_b$), $b$ edges connecting the previous sets to sets with $s-2$ elements (these edges have directions $e_{b+1}, e_{b+2}, \cdots, e_{2b}$), etc. After $n$ groups, we have used all $d = nb$ parallel directions; this process is repeated $t$ times as we cycle through all the directions and each group of $b$ directions connecting us to sets with one fewer element. In the last group of $b$, everything is connected to the empty set.

Although we have constructed $(T_S)_0$ starting from $[S]$ and working backward to $\emptyset$, for the next part of the proof, it is easier if we reorder the directions $\{e_i\}$ so that these disjoint paths can be broken into levels: the first level consisting of $b$ edges (those emanating from $\emptyset$) having directions $e_1, e_2, \cdots, e_b$; the second level consisting of edges having directions $e_{b+1}, e_{b+2}, \cdots, e_{2b}$; and so on. We call the $i^{\text{th}}$ level $S_i$.

We now have to label the graphs $(T_S)_x$ so that their union forms a communication graph. We do this by using Lemma 7.11 to associate to each equivalence class $[S]$ with $S$ a special subset, $n-1$ distinct equivalence classes $[\bar{S}_k]$ of regular $\bar{S}_k$. We now label $C_0 = (T_S)_0 \bigcup_{i=1}^{n-1} (T_{\bar{S}_k})_0$ as follows. We can think of each of the s levels of $(T_{\bar{S}_k})_0$ as being decomposed into $n$ groups of $b$ edges: the first group has directions $e_1, e_2, \cdots, e_b$; the second group has directions $e_{b+1}, e_{b+2}, \cdots, e_{2b}$ and so on. We call the $j^{\text{th}}$ group in the $i^{\text{th}}$ level $R_{ij}^k$. We have to assign $(n-1)s + t$ times to all these groups, $\{R_{ij}^k\}$ and $S_i$, in a manner which will satisfy the incomparability parts (i) and (ii) of the definition of task graph. This is fairly simple: we order the times so that they increase with increasing level. The only problem arises when we come to the "seam" between two levels. This scheme might assign the same time to an edge in one level as an edge in the next level which forms a directed path in one of the task graphs which makes up $C$. To remedy this, we make sure that the times on the seam are assigned to groups which come from $n$ distinct tasks; then we assign times to the remaining groups at each level. To make this more precise, let us reorganize the $n \times (n-1)$ groups $\{R_{ij}^k\}$ at each fixed level $i$. A typical collection of $n$ groups has the form $C_k = (R_{i1}^k, R_{i2}^{k+1}, \cdots, R_{in-1}^{k+n-2}, R_{in}^k)$ where the upper indices are read modulo $n-1$. Notice that the lower indices are all different - all directions are used exactly once in this collection; also, only the task $k$ appears twice. The set $\{R_{ij}^k\}$ is the disjoint union of the $n-1$ collections $C_k$; $S_i$ is the only other group which must be assigned a time at this level. We assign the times so that a collection $C_k$ is kept as contiguous as possible; every group in the collection is assigned one of two neighboring times. In addition, $S_i$ is inserted on the seam after all of the $\{R_{ij}^k\}$ and then the $\{R_{i+1 j}^k\}$ are cyclicly permuted so that they start with $R_{i+1 i}^{i+1}$ and continue. The details are messy to describe and are left to the reader; an example will probably be helpful. In this example, $n = 4$. The collections $C_k$ are $(R_{i1}^1, R_{i2}^2, R_{i3}^3, R_{i4}^1)$, $(R_{i1}^2, R_{i2}^3, R_{i3}^1, R_{i4}^2)$, $(R_{i1}^3, R_{i2}^1, R_{i3}^2, R_{i4}^3)$. When the $S_i$ are taken into account, we have the following collections on the seam: $(S_1, R_{22}^2, R_{23}^3, R_{24}^1)$, $(R_{21}^1, S_2, R_{33}^3, R_{34}^1)$, $(R_{31}^1, R_{32}^2, S_3, R_{44}^1)$. The full time list is given in Table 1 (after $n$ levels, we repeat the whole process; this happens $t$ times).

Let us summarize what we have done. We have made $C_0$ a communication graph with $(n-1)s + t$ times, which we can assume are consecutive odd integers starting with 1. Thus by the device of assigning the same times to $C_x$, when $x$ is even and assigning the following even times to $C_x$ when $x$ is odd, we can, by Lemma 7.4, assure that the union of $C_x$ for all $x \in Q_d$ is a communication graph with time $2((n-1)s + t)$. We note that this process requires $\dfrac{2((n-1)s + t)}{(n-1)d + b}$ times per element. Now

$$\frac{2((n-1)s+t)}{(n-1)d+b} = \frac{2(n(n-1)s+nt)}{n(n-1)d+nb} = 2\frac{(n(n-1)+1)s}{(n(n-1)+1)d} = 2\frac{s}{d}.$$

Thus the total time for all vertices at distance $s \leq d-1$ is $2\frac{s}{d}\binom{d}{s} = 2\binom{d-1}{s-1}$ as stated.

The algorithm for communicating simultaneously to vertices at distance $d-1$ and $d$ is a specialization of the method used above, but $n$ is set to 1. Details of the proof of this and the last statement of the theorem are left to the reader.

*Table I. Example with $n = 4$.*

| Time | List of Groups | |
|---|---|---|
| 1 | $R^1_{11}, R^2_{12}, R^3_{13}, R^1_{14}$ | |
| 2 | $R^2_{11}, R^3_{12}, R^1_{13}, R^2_{14}$ | |
| 3 | $R^3_{11}, R^1_{12}, R^2_{13}, R^3_{14}$ | |
| 4 | $S_1, R^2_{22}, R^3_{23}, R^1_{24}$ | Seam 1 |
| 5 | $R^2_{21}, R^3_{22}, R^1_{23}, R^2_{24}$ | |
| 6 | $R^3_{21}, R^1_{22}, R^2_{23}, R^3_{24}$ | |
| 7 | $R^1_{21}, S_2, R^3_{33}, R^1_{34}$ | Seam 2 |
| 8 | $R^2_{31}, R^3_{32}, R^1_{33}, R^2_{34}$ | |
| 9 | $R^3_{21}, R^1_{22}, R^2_{23}, R^3_{24}$ | |
| 10 | $R^1_{31}, R^2_{32}, S_3, R^1_{44}$ | Seam 3 |
| 11 | $R^2_{41}, R^3_{42}, R^1_{43}, R^2_{44}$ | |
| 12 | $R^3_{41}, R^1_{42}, R^2_{43}, R^3_{44}$ | |
| 13 | $R^1_{41}, R^2_{42}, R^3_{43}, S_4$ | |

Finally, we note that we have regularly ordered the hypercube.

*7.13. Theorem.* For two-way communication on the $d$-cube $Q^d$:

(i) $Q^d$ with the standard generators has a regular order;

(ii) the time for broadcast from each vertex to all vertices of distance $l$ is

$$\left\lceil \frac{\sum_{i=1}^{d}\binom{d}{i}-1}{d} \right\rceil;$$

(iii) the time for universal broadcast is
$$\left\lceil \frac{2^d - 1}{d} \right\rceil;$$

(iv) the time to communicate from each vertex to all vertices of distance $d-1$ and $d$ is $d$;

(v) the time to communicate from each vertex to all vertices of distance $s < d-1$ is
$$\binom{d-1}{s-1};$$

(vi) the time for a universal exchange is $2^{d-1}$.

*Proof.* This is exactly the content of Lemma 7.2, Corollary 7.5 and Lemma 7.7. Details are left to the reader (see the first part of the proofs of Theorems 6.3 and 6.5).

**§8. Summary.** We have discussed several combinatorial problems that arise when one looks at computational algorithms for highly symmetric networks of processors. We have developed a framework for dealing with these massively complicated scheduling tasks in a systematic way. More specifically, we have been concerned with minimal times associated with three communication tasks: universal broadcast, every processor has a vector that it wishes to broadcast to all the others; universal accumulation, every processor wishes to receive the sum of all the vectors being sent to it by all the other processors; and universal exchange, every processor wishes to exchange a vector with each other processor. We have exhibited some general properties that a network can have so that these tasks can be accomplished in optimal time.